\documentclass{amsart}
\usepackage{amssymb}
\usepackage{amsmath}
\usepackage{amscd}
\usepackage{url}
\usepackage[english]{babel}
\newtheorem{prop}{Proposition}[section]

\newtheorem{conj}[prop]{Conjecture}

\theoremstyle{definition}

\numberwithin{equation}{section}

\renewcommand{\Im}{{\mathrm {Im}}}
\newcommand{\ord}{{\mathrm {ord}}}

\newcommand{\Hom}{{\mathrm {Hom}}}

\newcommand{\Fil}{{\mathrm {Fil}}}

\newcommand{\alg}{\mathrm{alg}}

\newcommand{\Tr}{{\mathrm {Tr}}}

\newcommand{\Sym}{{\mathrm {Sym}}}

\newcommand{\Spin}{{\mathrm {Spin}}}
\newcommand{\dR}{{\mathrm {dR}}}

\newcommand{\Frob}{{\mathrm {Frob}}}

\newcommand{\Gal}{\mathrm {Gal}}

\newcommand{\A}{{\mathbb A}}
\newcommand{\CC}{{\mathbb C}}
\newcommand{\C}{{\mathbb C}}
\newcommand{\RR}{{\mathbb R}}

\newcommand{\QQ}{{\mathbb Q}}
\newcommand{\Q}{{\mathbb Q}}
\newcommand{\ZZ}{{\mathbb Z}}

\newcommand{\HH}{{\mathfrak H}}

\newcommand{\VVV}{{\mathbb V}}

\newcommand{\n}{{\mathfrak n}}

\newcommand{\q}{{\mathfrak q}}

\newcommand{\FF}{{\mathbb F}}

\newcommand{\GL}{\mathrm {GL}}
\newcommand{\PGL}{\mathrm {PGL}}

\newcommand{\SL}{\mathrm {SL}}
\newcommand{\Sp}{\mathrm {Sp}}
\newcommand{\SO}{\mathrm {SO}}
\newcommand{\GSp}{\mathrm {GSp}}
\newcommand{\PGSp}{\mathrm {PGSp}}

\newcommand{\Qbar}{\overline{\mathbb Q}}

\newcommand{\rhobar}{\overline{\rho}}

\begin{document}
\title{$\GL_2\times\GSp_2$ $L$-values and Hecke eigenvalue congruences}
\author{Jonas Bergstr\"om}
\address{Matematiska institutionen\\ Stockholms universitet\\ 106 91 Stockholm\\Sweden.}
\email{jonasb@math.su.se}

\author{Neil Dummigan}
\address{University of Sheffield\\ School of Mathematics and Statistics\\
Hicks Building\\ Hounsfield Road\\ Sheffield, S3 7RH\\
U.K.}
\email{n.p.dummigan@shef.ac.uk}

\author{David Farmer}
\address{American Institute of Mathematics\\600 East Brokaw Road\\San Jose, CA 95112\\U.S.A.}
\email{farmer@aimath.org}

\author{Sally Koutsoliotas}
\address{Department of Physics and Astronomy\\Bucknell University\\Lewisburg, PA 17837\\U.S.A.}
\email{koutslts@bucknell.edu}


\begin{abstract}
\hspace{-4.5pt}
We find experimental examples of congruences of Hecke eigenvalues between automorphic representations of groups such as
$\GSp_2(\A)$, $\SO(4,3)(\A)$ and $\SO(5,4)(\A)$, where the prime modulus should, for various reasons, appear in the algebraic part of a critical ``tensor-product'' $L$-value associated to cuspidal automorphic representations of $\GL_2(\A)$ and $\GSp_2(\A)$. Using special techniques for evaluating $L$-functions with few known coefficients, we compute sufficiently good approximations to detect the anticipated prime divisors.
\end{abstract}

\maketitle

\section{Introduction}
This paper is about divisors of critical values of $L$-functions as moduli of congruences between Hecke eigenvalues of automorphic forms. It is made possible by three separate developments in computational number theory.
\begin{enumerate}
\item M\'egarban\'e's large-scale computation of traces of Hecke operators on spaces of level-one algebraic modular forms, for $\SO(7),\,\SO(8)$ and $\SO(9)$ \cite{Me}, following the endoscopic classification of the associated automorphic representations by Chenevier and Renard \cite{CR}.
\item Faber and van der Geer's computation of traces of Hecke operators on spaces of vector-valued Siegel modular forms of genus $2$ and level one, using point counts on hyperelliptic curves of genus $2$ over finite fields \cite{FvdG}. When the space is $1$-dimensional, this gives Hecke eigenvalues. They computed traces of Hecke operators $T(p)$ and $T(p^2)$ for $p,p^2\leq 37$. The first-named author of this paper refined their method and extended the bound to $179$. The data is available at \url{http://smf.compositio.nl/}.
\item A new technique for computing good approximations to values of $L$-functions satisfying functional equations, given only a few coefficients in the Dirichlet series, developed by the third-named author and Ryan \cite{FR}. This combines,  in such a way as to make unknown errors cancel, approximations obtained by the method of Rubinstein \cite{Rub}, which is related to the technique implemented in Magma \cite{Mag}, which is described in \cite{Do}.
\end{enumerate}
Conjecture 4.2 of \cite{BD} is a very wide generalisation of Ramanujan's mod $691$ congruence, to ``Eisenstein'' congruences between Hecke eigenvalues of automorphic representations of $G(\A)$, where $\A$ is the adele ring and $G/\Q$ is any connected, split reductive group. On one side of the congruence is a cuspidal automorphic representation $\tilde{\Pi}$. On the other is one induced from a cuspidal automorphic representation $\Pi$ of the Levi subgroup $M$ of a maximal parabolic subgroup $P$.
The modulus of the congruence comes from a critical value of a certain $L$-function, associated to $\Pi$ and to the adjoint representation of the $L$-group $\hat{M}$ on the Lie algebra $\hat{\n}$ of the unipotent radical of the maximal parabolic subgroup $\hat{P}$ of $\hat{G}$. Starting from $\Pi$, we conjecture the existence of $\tilde{\Pi}$, satisfying the congruence.

In \cite{BDM} we already used M\'egarban\'e's data for $\SO(7)$ and $\SO(8)$, observing experimental Eisenstein congruences for the cases $G=\SO(4,3)$, $M=\GL_1\times\SO(3,2)$ and $G=\SO(4,4),\,M=\GL_2\times\SO(2,2)$. To support the conjecture, we then needed to find the observed moduli in the corresponding $L$-values. In the $\SO(4,4)$ case these were triple product $L$-values for elliptic modular forms, which were computed exactly by Ibukiyama and Katsurada, using the pull-back of a genus $3$ Siegel-Eisenstein series. In the $\SO(4,3)$ case they were spinor $L$-values for vector-valued Siegel modular forms of genus $2$ and level one (note that $\SO(3,2)\simeq \PGSp_2$), and we resorted to sufficiently good numerical approximations. For this, Magma was good enough, and we used the Hecke eigenvalues computed by the first-named author, which went as far as the bound $149$ at that point.

For this paper we used M\'egarban\'e's data for $\SO(9)$ to find an experimental Eisenstein congruence (mod $q=17$) in the case $G=\SO(5,4),\,M=\GL_2\times\SO(3,2)$ (see Example 5 in \S 7).
In this case, the associated $L$-function has degree-$8$ Euler factors, and is the ``tensor-product'' of the Hecke $L$-function of an elliptic modular form and the spinor $L$-function of a vector-valued Siegel modular form of genus $2$ (both level one). We also used M\'egarban\'e's data for $\SO(7)$ to find experimental ``endoscopic'' congruences (mod $q=71$ and $61$), between functorial lifts from $\SO(2,1)\times\SO(3,2)$ to $\SO(4,3)$ and non-lifts on $\SO(4,3)$ (see Examples 3 and 4 in \S 6, and note that $\SO(2,1)\simeq\PGL_2$). Here, as with the Eisenstein congruences, a generalisation of a construction of Ribet leads from the congruence, via an extension of mod $q$ Galois representations, to an element of order $q$ in a Bloch-Kato Selmer group, which according to the Bloch-Kato conjecture ought to show up in a certain critical $L$-value. For these $\SO(7)$ endoscopic congruences, again it is a tensor-product $\GL_2\times\GSp_2$ $L$-function, for an elliptic modular form and a vector-valued Siegel modular form of genus $2$. These congruences are analogous to those between Yoshida lifts and non-lifts (Siegel modular forms of genus $2$) appearing in \cite{BDS}, where the $L$-function is a degree-$4$ tensor-product $L$-function for two elliptic modular forms.

To obtain sufficiently good approximations to the $\GL_2\times\GSp_2$ $L$-values, Magma requires many more coefficients in the Dirichlet series than we could obtain using the computations of Hecke eigenvalues for Siegel modular forms. (See \cite[\S 7]{BDM} for a comment on the difficulty of extending these much further.) So for this the third and fourth-named authors used the kind of averaging technique described in \cite{FR}. This is described in \S 4.3--\S 4.6, and the results are in \S 4.2. The numerical approximations to ratios of $L$-values (and appropriate powers of $\pi$) are very close to rational numbers, and in these rational numbers we find the expected factors of $17,71$ and $61$. We also stumbled on some other factors of $839$ and $61$ (again), and realised at that point that these could also be explained using the Bloch-Kato conjecture, in terms of Eisenstein congruences for $G=\GSp_4$. For $q=839$ (Example 1), $P$ is the Klingen parabolic subgroup (Kurokawa-Mizumoto congruences) while for $q=61$ (Example 2), $P$ is the Siegel parabolic subgroup (Harder's conjecture). This is explained in Section 5. In fact, following basic background in \S 2, we begin with these examples in \S 3, where a rougher heuristic is given. As noted in \S \ref{KMt}, \S\ref{Hat}, the accidental discovery of these experimental divisibilities led to the {\em proof } of the analogous divisibilities in the scalar-valued case, where the necessary pullback formulas are known. However, the experimental congruences and divisibilities of $L$-values involving $\SO(7)$ and $\SO(9)$ are presumably some way beyond what can currently be proved (with the exception of the congruence in \S 6, Example 4). This justifies the effort made to take the experimental results as far as possible, extending to types of congruences and numerical techniques not considered or employed in \cite{BDM}.

In \S 9 we use M\'egarban\'e's $\SO(9)$ data to observe an experimental endoscopic congruence mod $q=37$, between a functorial lift from $\SO(3,2)\times\SO(3,2)$ to $\SO(5,4)$ and a non-lift on $\SO(5,4)$. Applying the same techniques to the associated $\GSp_2\times\GSp_2$ $L$-values, we cannot obtain such reliable and accurate approximations as before, but are still able to see the expected factor of $37$.

As should already be clear from this introduction, the data that was computed by Thomas M\'egarban\'e was indispensable, and we are grateful to him for providing it to us before it was publicly available. We thank also Mark Watkins for pointing out errors in an earlier version of this paper, including one that had a very significant impact on \S \ref{SO9}, and for corroborating the new computation thus enabled.

\section{$\GL_2\times\GSp_2$ $L$-functions}
Let $f$ be a normalised cuspidal Hecke eigenform of weight $\ell$ for $\SL_2(\ZZ)$. Then $f\left(\frac{a\tau+b}{c\tau+d}\right)=(c\tau+d)^{\ell}f(\tau)$ for all $\begin{pmatrix}a & b\\c & d\end{pmatrix}\in\SL_2(\ZZ)$, and $\tau\in\HH=\{\tau\in\CC:\,\Im(\tau)>0\}$, and $f(\tau)=\sum_{n=1}^{\infty}a_n(f)q^n$, with $q=e^{2\pi i\tau}$ and $a_1=1$. The Fourier coefficients are also the eigenvalues of Hecke operators. The Hecke $L$-function is $$L(s,f)=\prod_{p\,\text{prime}}(1-a_p(f)p^{-s}+p^{\ell-1-2s})^{-1}.$$
Let $1-a_p(f)X+p^{\ell-1}X^2=:(1-\alpha_{p,1}X)(1-\alpha_{p,2}X)$.

Let $F$ be a cuspidal Hecke eigenform of weight $\Sym^j\otimes\det^k$ for $\Sp_2(\ZZ):=\{g\in M_4(\ZZ):\,{}^tgJg=J\}$, where $J=\begin{pmatrix} 0_2 & -I_2\\I_2 & 0_2\end{pmatrix}$. Then $F:\HH_2\rightarrow V$, where $\HH_2=\{Z\in M_2(\CC):{}^tZ=Z, \Im(Z)>0\}$ is Siegel's upper half space of genus $2$, $V$ is the space of the representation $\rho=\Sym^j(\CC^2)\otimes\det^k$ of $\GL_2(\CC)$, and $$F\left((AZ+B)(CZ+D)^{-1}\right)=\rho(CZ+D)(F(Z))\,\,\,\forall \begin{pmatrix} A & B\\C & D\end{pmatrix}\in \Sp_2(\ZZ).$$
Let the elements $T(p), T(p^2)$ of the genus-$2$ Hecke algebra be as in \cite[\S 16]{vdG} (with the scaling as following Definition 8). Let $\lambda_F(p), \lambda_F(p^2)$ be the respective eigenvalues for these operators acting on $F$.
The spinor $L$-function of $F$ is $L(s,F,\Spin)=\prod_{p\,\text{prime}}L_p(s,F,\Spin)$,
where
\begin{multline*} L_p(s,F,\Spin)^{-1}=1-\lambda_F(p)p^{-s}+(\lambda_F(p)^2-\lambda_F(p^2)-p^{j+2k-4})p^{-2s}  \\
-\lambda_F(p)p^{j+2k-3-3s}+p^{2j+4k-6-4s}.
\end{multline*}
Let $L_p(s,F,\Spin)^{-1}=:\prod_{j=1}^4(1-\beta_{p,j}p^{-s})$.

Now we define $L(s,f\otimes F):=\prod_{p\,\text{prime}}L_p(s,f\otimes F),$
where $$L_p(s,f\otimes F)^{-1}:=\prod_{i=1}^2\prod_{j=1}^4(1-\alpha_{p,i}\beta_{p,j}p^{-s}).$$ To understand the conjectured functional equation and critical values for this $L$-function, it is convenient to introduce the motive $M_f$ attached to $f$, and the conjectured motive $M_F$ attached to $F$, of ranks $2$ and $4$ respectively. The Betti realisations have Hodge decompositions $M_{f,B}\otimes\CC\simeq H^{0,\ell-1}\oplus H^{\ell-1,0}$ and $M_{F,B}\otimes\CC\simeq H^{0,j+2k-3}\oplus H^{j+2k-3,0}\oplus H^{k-2,j+k-1}\oplus H^{j+k-1,k-2}$, with each $H^{p,q}$ $1$-dimensional. The $L$-functions associated to ($q$-adic realisations of) $M_f$ and $M_F$ are $L(s,f)$ and $L(s,F,\Spin)$ respectively. The $L$-function $L(s,f\otimes F)$ is associated to the rank-$8$ motive $M:=M_f\otimes M_F$, which has Hodge decomposition $M_B\otimes\CC\simeq\oplus(H^{p,q}\oplus H^{q,p})$, where $p+q=j+2k+\ell-4$ and $p\in\{0,k-2,\min\{\ell-1,j+2k-3\},\min\{j+k-1,k+\ell-3\}\}=:\{p_1,p_2,p_3,p_4\}$, where we label the elements so that $p_1\leq p_2\leq p_3\leq p_4$. According to \cite[Table 5.3]{De}, each $(p,q)$ contributes $i^{q-p+1}$ to the sign in the conjectural functional equation, and using the fact that $j$ is even, one checks easily that the sign should be $+1$. Following the recipe in \cite{Se} (or see again \cite[Table 5.3]{De}), the product of gamma factors is
$\gamma(s)=\prod_{i=1}^4\Gamma_{\CC}(s-p_i)$, where $\Gamma_{\CC}(s):=(2\pi)^{-s}\Gamma(s)$. Note that, following \cite[Remark 6.2]{BH1}, it makes no difference to replace any $p_i$ by $q_i=j+2k+\ell-4-p_i$. Anyway, the conjectured functional equation is $\Lambda(j+2k+\ell-3-s)=\Lambda(s),$ where $\Lambda(s):=\gamma(s)L(s,f\otimes F)$. The meromorphic continuation and functional equation have been proved by B\"ocherer and Heim \cite{BH1} in the case that $F$ is scalar valued (i.e. $j=0$), Furusawa \cite{Fu} having already dealt with the even more special case $\ell=k$ (and $j=0$).

The critical values are $L(t,f\otimes F)$ for integers $t$ such that neither $\gamma(s)$ nor $\gamma(j+2k+\ell-3-s)$ has a pole at $s=t$. This is for $p_4<t\leq q_4$. In all our examples, the coefficient field of $M_f$ and $M_F$ (hence of $M$) is $\QQ$, so we suppose for convenience that this is the case. (Then $M_B$ and $M_{\dR}$ are $\QQ$-vector spaces.) For each critical $t$, there is a Deligne period $c^+(M(t))$ defined as in \cite{De}, up to $\QQ^{\times}$ multiples. (It is the determinant, with respect to bases of $4$-dimensional $\Q$-vector spaces $M_B(t)^+$ and $M_{\dR}(t)/\Fil^0$, of an isomorphism between $M_B(t)^+\otimes\CC$ and $(M_{\dR}(t)/\Fil^0)\otimes\CC$.) Deligne's conjecture (in this instance) is that $L(s,f\otimes F)/c^+(M(t))\in \Q^{\times}$. Later we shall sometimes make a special choice of $c^+(M(t))$, and define $L_{\alg}(t,f\otimes F)=L(t,f\otimes F)/c^+(M(t))$. If $t,t'$ are critical points with $t\equiv t'\pmod{2}$, then $c^+(M(t'))=(2\pi i)^{4(t'-t)}c^+(M(t))$, because $M_B(t')=M_B(t)(2\pi i)^{t'-t}$ while $M_{\dR}(t)/\Fil^0$ does not change for $t$ within the critical range. So the ratio $\frac{L_{\alg}(t',f\otimes F)}{L_{\alg}(t,f\otimes F)}=\frac{L(t',f\otimes F)}{(2\pi i)^{4(t'-t)}L(t,f\otimes F)}$, which should be a rational number, is independent of any choices. Prime divisors of its numerator or denominator will turn out to be significant. Actually, for these tensor product motives the condition $t\equiv t'\pmod{2}$ is unnecessary \cite[Cor.1, p.1188]{Y}, but among our examples, only for the one in \S \ref{SO9} are there so few critical points that we would need to relax it.

\section{Expected consequences of congruences: the rough version}
\subsection{Kurokawa-Mizumoto type}\label{KMt}
Suppose that $\ell=j+k$, with $\ell$ as above, $j,k$ non-negative even integers. There is a vector-valued Klingen-Eisenstein series $[f]_j$, a non-cuspidal genus-$2$ Siegel modular form of weight $\Sym^j\otimes\det^k$ for $\Sp_2(\ZZ)$, satisfying $\Phi([f]_j)=f$, where $\Phi$ is the Siegel operator. (See \cite[\S 1]{A} for more details.) Let $q>2\ell$ be a prime divisor of the numerator of $L_{\alg}(2\ell-2-j,\Sym^2 f)$, which we can take to be $L(2\ell-2-j,\Sym^2 f)/\pi^{2k+\ell-3}(f,f)$, where $(f,f)$ is the Petersson norm and $L(s,\Sym^2 f)=\prod_{p\,\text{prime}}((1-\alpha_1^2p^{-s})(1-\alpha_1\alpha_2p^{-s})(1-\alpha_2^2p^{-s}))^{-1}$. Sometimes it is possible to prove a congruence (mod $q$) of Hecke eigenvalues between $[f]_j$ and some cuspidal Hecke eigenform $F$, also of weight $\Sym^j\otimes\det^k$ for $\Sp_2(\ZZ)$. The first examples were proved by Kurokawa and Mizumoto \cite{K,Mi}, with further examples proved by Satoh \cite{Sa} and in \cite{Du2}.

Note that on $[f]_j$ the eigenvalue of $T(p)$ is $a_p(f)(1+p^{k-2})$, in fact its spinor $L$-function (defined in terms of Hecke eigenvalues just as for the cuspidal case) is
$L(s,[f]_j,\Spin)=L(s,f)L(s-(k-2),f)$. Then $L(s,f\otimes [f]_j)=L(s,f\otimes f)L(s-(k-2),f\otimes f)$. Since $L(s,f\otimes f)=\zeta(s-(\ell-1))L(s,\Sym^2 f)$, we find that $$L(s,f\otimes [f]_j)=\zeta(s-(\ell-1))\zeta(s-(\ell+k-3))L(s,\Sym^2 f)L(s-(k-2),\Sym^2 f).$$
In this situation where $\ell=j+k$, the critical range for $L(s,f\otimes F)$ is $\ell\leq t\leq\ell+k-3$.
The factor $\zeta(s-(\ell+k-3))$ is non-zero at $s=\ell+k-3$, but has a trivial zero at all other odd $s$ in the critical range, e.g. at $s=\ell+k-5$. Checking the other factors, we find that $\frac{L(s,f\otimes[f]_j)}{\pi^8L(s-2,f\otimes[f]_j)}$ has a simple pole at $s=\ell+k-3$.

The mod $q$ congruence of Hecke eigenvalues between $[f]_j$ and $F$, hence between coefficients of the Dirichlet series for $L(s,f\otimes[f]_j)$ and $L(s,f\otimes F)$, might lead one roughly to expect that the pole of $\frac{L(s,f\otimes[f]_j)}{\pi^8L(s-2,f\otimes[f]_j)}$ at the rightmost critical point $s=\ell+k-3$ (for $L(s,f\otimes F)$) should cause a pole {\em mod $q$} of $\frac{L(j+2k-3,f\otimes F)}{\pi^8L(j+2k-5,f\otimes F)}$, i.e. a factor of $q$ in its denominator. Since $j+2k-5$ could be replaced by any odd $s$ in the critical range strictly to the left of $j+2k-3$, we can think of this $q$ as being in the denominator of $L_{\alg}(j+2k-3,f\otimes F)$, without worrying too much about the correct scaling.

We have so far chosen the ``motivic'' normalisation of the $L$-function, but it is also convenient to consider the ``unitary'' normalisation
$L(s+(j+2k+\ell-4)/2,f\otimes F)$, which should satisfy a functional equation relating $s$ and $1-s$. This normalisation is natural if we consider the $L$-function as arising from automorphic representations $\pi_f$ and $\pi_F$ of $\GL_2(\A)$ and $\GSp_2(\A)$ respectively, where $\A$ is the adele ring, so we set $L(s,\pi_f\otimes\pi_F):=L(s+(j+2k+\ell-4)/2,f\otimes F)$. We then expect to find $q$ in the denominator of $\frac{L((k-2)/2,\pi_f\otimes \pi_F)}{\pi^8L((k-6)/2,\pi_f\otimes \pi_F)}$.

{\bf Example 1. } $\ell=16, j=4, k=12, q=839$. The congruence is \cite[Proposition 4.1]{Du2}. We expect $\frac{L(5,\pi_f\otimes \pi_F)}{\pi^8L(3,\pi_f\otimes \pi_F)}$ to be a rational number with $839$ in the denominator. In fact, if we observe it in the denominator of
$\frac{L(5,\pi_f\otimes \pi_F)}{\pi^8L(3,\pi_f\otimes \pi_F)}$ (or equivalently in the numerator of $\frac{\pi^8L(3,\pi_f\otimes \pi_F)}{L(5,\pi_f\otimes \pi_F)}$), but not in the denominator of $\frac{\pi^8L(1,\pi_f\otimes \pi_F)}{L(3,\pi_f\otimes \pi_F)}$, we should feel reasonably confident that it is coming from the denominator of $L_{\alg}(5,\pi_f\otimes\pi_F)$ rather than the numerator of $L_{\alg}(3,\pi_f\otimes\pi_F)$, since it seems unlikely that it would also happen to divide the numerator of $L_{\alg}(1,\pi_f\otimes\pi_F)$.

In our numerical examples in \S 4, $j>0$ and $F$ is vector-valued. However, in the case that $j=0$ and $F$ is scalar-valued, a formula of Heim, for the restriction of a genus $5$ Eisenstein series to $\HH_1\times\HH_2\times\HH_2$ \cite{He}, in which appears the $L$-value in question, allows one to actually prove the expected divisibility \cite{DHR}.
\subsection{Harder type}\label{Hat} Suppose that $k'=j+2k-2$, with $k'$ the weight of a normalised, cuspidal Hecke eigenform $g$ for $\SL_2(\ZZ)$, $j> 0, k\geq 3$ integers with $j$ even. In many examples there appears to be a cuspidal Hecke eigenform $F$, of weight $\Sym^j\otimes\det^k$ for $\Sp_2(\ZZ)$, and a congruence
$\lambda_F(p)\equiv a_p(g)+p^{k-2}+p^{j+k-1}\pmod{\q}$, where $\q$ divides the numerator of a suitably normalised $L_{\alg}(j+k,g)$. Cases where such congruences have been verified for $p\leq 37$, using Hecke eigenvalue computations by Faber and van der Geer \cite{FvdG}, are described in \cite{vdG}. The original example $(k',j,k,q)=(22,4,10,41)$ used by Harder to support his conjecture \cite{H}, has subsequently been proved by Chenevier and Lannes \cite[Chapter 10, Theorem* 4.4(1)]{CL}.

One way of expressing the congruence is to say that $L(s,F,\Spin)$ is congruent, coefficient by coefficient, to $L(s,g)\zeta(s-(k-2))\zeta(s-(j+k-1))$. Then, with auxiliary $f$ of weight $\ell$, $L(s,f\otimes F)$ is congruent, coefficient by coefficient, to $L(s,f\otimes g)L(s-(k-2),f)L(s-(j+k-1),f)$. Now if $\ell/2$ is odd then the sign in the functional equation of $L(s,f)$ is $-1$, so $L(\ell/2,f)=0$ and the factor $L(s-(j+k-1),f)$ vanishes at $s=(\ell/2)+j+k-1$. We might then roughly expect $L_{\alg}((\ell/2)+j+k-1,f\otimes F)$ to vanish {\em mod $q$}, so to find $q$ in the numerator of $\frac{L((\ell/2)+j+k-1,f\otimes F)}{\pi^8L((\ell/2)+j+k-3,f\otimes F)}$, which is the same as $\frac{L((j+2)/2,\pi_f\otimes\pi_F)}{\pi^8L((j-2)/2,\pi_f\otimes\pi_F)}$. These will be critical values as long as we ensure that $\ell\geq j+4$.

{\bf Example 2. } If $(j,k)=(4,15)$, so that $k'=32$, then the unique (up to scaling) cusp form $F$ of weight $\Sym^4\otimes\det^{15}$ appears to satisfy a congruence as above mod $\q$, with $\q\mid L_{\alg}(19,g)$ a divisor of $61$.
(Note that there are two conjugate choices for $g$ and for $\q$.) Now let $f$ be the unique normalised cusp form of weight $\ell=18$ for $\SL_2(\ZZ)$. We expect to find $61$ in the numerator of $\frac{L(3,\pi_f\otimes\pi_F)}{\pi^8L(1,\pi_f\otimes\pi_F)}$.

In our numerical examples in \S 4, $j>0$ and $F$ is vector-valued. However, in the case that $j=0$ and $F$ is scalar-valued, congruent mod $\q$ to the Saito-Kurokawa lift of $g$, a formula of Saha, for the restriction to (Siegel) $\HH_1\times\HH_2$ of a genus $3$ Hermitian Eisenstein series (non-holomorphic and non-convergent in our case), allows one to actually prove the expected divisibility. This is in the recent Sheffield Ph.D. thesis of Rendina \cite{Re}.

\section{Computing the $L$-values}
\subsection{Generalities}\label{computingAp}
The method for computing Hecke eigenvalues of genus $2$ cusp forms $F$, hence coefficients of $L(s,F,\Spin)$,
is described in \cite[\S 7]{BDM}.
It was first carried out for $p\leq 37$ by Faber and van der Geer \cite{FvdG},
and extended by the first-named author to obtain the first $180$ coefficients in the Dirichlet series.
As described in those references, computing the $p$th Dirichlet coefficient
requires approximately $p^4$ operations.  Thus, computing significantly more
Dirichlet coefficients is not practical.

\subsection{Computational results}\label{compres}
Using the method described in Section~\ref{sec:howtoevaluateL},
we have experimentally determined the following expressions involving special values.

\vskip 0.1in
\noindent Case 1: $\ell=\mathbf{16}, (j,k)=\mathbf{(4,12)}$
\vskip 0.1in

We have
$$
\frac{\pi^8L(3,\pi_f\otimes \pi_F)}{L(5,\pi_f\otimes \pi_F)}
=
\frac{7^2 \cdot 17 \cdot 839}{2^3 \cdot 3^2}
$$
and
$$
\frac{\pi^8L(1,\pi_f\otimes \pi_F)}{L(3,\pi_f\otimes \pi_F)}
=
\frac{3^4 \cdot 7 \cdot 11^2 \cdot 71}{2^2 \cdot 5 \cdot 17}.
$$
The $839$ is as predicted by Example 1 in \S \ref{KMt}. The $71$ will be explained in Example~3 in \S \ref{endo7} below, as will the $17$ in Example 5 in \S \ref{eis9} below.

\vskip 0.1in
\noindent Case 2: $\ell=\mathbf{18}, (j,k)=\mathbf{(4,15)}$

We have
$$
 \frac
{\pi^8L(1,\pi_f\otimes\pi_F)}
{L(3,\pi_f\otimes\pi_F)}
=
\frac
{3 \cdot 5^4 \cdot 7 \cdot 13^2 \cdot 193}
{2^4 \cdot 11 \cdot 61}.
$$
The $61$ is as predicted by Example 2 in \S \ref{Hat}. For a comment on why we are unable to account for the $193$, see the end of Example 3 in \S \ref{endo7}.

\vskip 0.1in
\noindent Case 3:  $\ell=\mathbf{16}, (j,k)=\mathbf{(6,10)}$

We have
$$
 \frac{\pi^8L(1,\pi_f\otimes \pi_F)}{L(3,\pi_f\otimes \pi_F)}
=
\frac{3^4 \cdot 5^2 \cdot 61}{2^3}.
$$
The $61$ in the numerator will be revisited in Example 4 in \S \ref{endo7} below.

\subsection{Numerically evaluating L-functions}\label{sec:howtoevaluateL}

We describe the numerical evaluation of the degree 8 L-function
$L(s, \pi_f\otimes\pi_F)$, using the case
$\ell=16$, $(j,k)=(6,10)$ as a representative example.
We wish to evaluate that L-function at $s=1$ and $s=3$ to
high precision so as to confidently identify
a normalized ratio of those values as a rational
number.  (The meaning of `high precision' depends on context.
Here we will consider 30 digits to be
a reasonable target.)

First we explain why this requires some effort.  It is straightforward
to make as many Dirichlet coefficients of $L(s, \pi_f)$ as we wish,
but for $L(s, \pi_F)$ we have only the first 180 coefficients,
so we have only the first 180 coefficients of
$L(s) := L(s, \pi_f\otimes\pi_F)$.
As remarked in Section~\ref{computingAp}, it is prohibitive to produce
significantly more coefficients.

A common method of evaluating L-functions is to use the built-in
functionality of Magma~\cite{Mag}.  Since the functional equation satisfied
by this L-function, in the unitary/analytic normalization, is
\begin{equation}
\Lambda(s) := \Gamma_\C(s+19)\Gamma_\C(s+11)\Gamma_\C(s+4)^2 L(s)
=
   \Lambda(1-s),
\end{equation}
Magma can tell us how many coefficients are needed:
\begin{verbatim}
> L:=LSeries(1,[4,4,5,5,11,12,19,20],1,0: Sign:=1);
> N:=LCfRequired(L);N;
4145
\end{verbatim}
Thus, more than 4000 coefficients are required for evaluating the L-function
using standard methods, but we only have~180.
Instead, we will use the methods of~\cite{FR}
to accurately evaluate the critical values using only the available coefficients.
We summarize the ideas as applied to this example.

High-precision evaluations of general L-functions use the so-called
approximate functional equation (see~\cite{Rub} for details
and technical conditions).  If $L(s)=\sum b_n n^{-s}$
has an analytic continuation to an entire function that satisfies the
functional equation $\Lambda(s) = G(s)L(s) = \varepsilon \Lambda(1-s)$,
and $g(s)$ is a suitable auxiliary function, then
\begin{equation}
\label{eqn:appfe}
  L(s) =
\sum_{n=1}^{\infty} \left(h_1(s,n) + \varepsilon h_2(1-s,n)\right) b_n
\end{equation}
where
\begin{align}\label{eqn:mellin}
   h_1(s,n) &:= (g(s)G(s))^{-1} \frac{1}{2\pi i}
          \int_{\nu - i \infty}^{\nu + i \infty}
         \frac{g(s+z)G(s+z)}{n^{s+z}} \frac{dz}{z}
                    \notag \\
  h_2(1-s,n) &:= (g(s)G(s))^{-1} \frac{1}{2\pi i}
                   \int_{\nu - i \infty}^{\nu + i \infty}
                     \frac{g(1-s+z)G(1-s+z)}{n^{1-s+z}} \frac{dz}{z} .
\end{align}
Here $\nu$ is any real number to the right of all poles of the integrand,
and $g(s)$ is any entire function such that the integrals converge absolutely.

For this discussion,
the important parameter is the test function
$g(s)$.
The idea is that we can evaluate \eqref{eqn:appfe} multiple times
with different test functions.  Each evaluation provides
slightly different information, which we can combine to overcome the
fact that we only have a few Dirichlet coefficients.
We will use test functions of the form $g(s) = e^{i\beta s + \alpha s^2}$ with
$\alpha > 0$, or $\alpha = 0$ and
$|\beta| < d \frac{\pi}{4}$, where $d$ is the degree of the L-function
(in our example, $d=8$).

If we insert a particular value for $s$, say $s=1$, and let
$g(s) = e^{i\beta s + s^2/1000}$, then~\eqref{eqn:appfe} has the form
\begin{equation}\label{eqn:genbeta}
L_\beta(1) = \sum_{n=1}^\infty c_\beta(n) b_n ,
\end{equation}
where $c_\beta(n)$ is a number which depends on $\beta$, $n$,
and the parameters in the functional equation.
The $L_\beta(1)$ on the left side of \eqref{eqn:genbeta} is independent of $\beta$,
but we use that notation to indicate which auxiliary function was used on the right side.
We can evaluate
$c_\beta(n)$ as accurately as we wish by numerically evaluating the
integrals that appear in~\eqref{eqn:mellin}, which we now describe.

\subsection{Numerically evaluating the integrals}
This section is a summary of material from~\cite{Rub}.
Our goal here is to provide information for someone
to reproduce our calculation; our goal is not to provide a detailed
exposition on numerically approximating integrals.

We wish to numerically calculate, to high precision, the numbers $c_\beta(n)$ in \eqref{eqn:genbeta}.
This involves evaluating an integral of a product of $\Gamma$-functions,
exponentials, and powers.  This can be done using any computer algebra package
which contains those functions and which can work to arbitrary precision;
our computations were done in Mathematica \cite{Mat}.

The main issues involved are:
\begin{enumerate}
\item \label{issue:precision} Evaluating the integrand to sufficiently high precision,
\item \label{issue:truncate} Truncating the improper integral, and
\item \label{issue:sum} Evaluating the resulting finite integral.
\end{enumerate}
For issue (\ref{issue:precision}) the concern is possible loss of significant digits due to
catastrophic cancellation.  This is a minor issue in the case at hand, although
it becomes a serious problem when evaluating L-functions at $s=\sigma + i t$
with $t$ large (because the completed L-function $\Lambda(s)$ is what is actually
being computed, and it decreases exponentially as a function of $t$).
See Section~3.3 of~\cite{Rub} for details.

In a system such as Mathematica, issue~(\ref{issue:precision}) is even less of a
problem, because the software keeps track of the precision of the calculation.
Should there be insufficient precision at the end, one merely re-calculates,
beginning with a larger precision.

As explained in Sections~2.4 and~3.6
of~\cite{Rub}, these integrals can be evaluated by a simple Riemann sum
(i.e., a sum that evaluates the integrand at equally spaced points).
In particular, the step size is inversely proportional to the number
of digits of accuracy in the result.
That addresses issue~(\ref{issue:sum}), and since the integrand
is decreasing exponentially, also issue (\ref{issue:truncate}).

\subsection{A numerical example}
We will evaluate \eqref{eqn:genbeta} with $\beta=0$.  We set $\nu=3$,
a stepsize of $1/5$ in the Riemann sum, summing from -29 to 29,
and evaluate the integrand
to 40 digits of precision.  We find:
\begin{align}\label{eqn:beta0}
L_0(1)  =  &\mathstrut  1.245\, b_1 + 0.534\, b_2 + 0.269\, b_3 + \cdots +
0.000668\, b_{17} \notag \\
&+ \cdots + 2.10\times 10^{-10}\, b_{101}  +
\cdots +
8.56\times 10^{-14}\, b_{181} \notag \\
& +
\cdots
+1.1\times 10^{-21}\, b_{499}
+
\cdots
+
5.5\times 10^{-29}\, b_{1009} \notag  \\
& +
\cdots
+
7.3 \times 10^{-34}\, b_{1499}
+
\cdots
+
8.3 \times 10^{-38}\, b_{1999} \notag  \\
&+
\cdots + 6.7 \times 10^{-53}\, b_{4999}  + \cdots.
\end{align}

One sees immediately that using 180 terms gives an error of
more than $10^{-13}$, which is far from our goal of 30 decimal digits.
Note that the numerical values in~\eqref{eqn:beta0}, and in all following
equations, are truncations of the actual value.  For example,
the actual computed
coefficient of $b_1$ in \eqref{eqn:beta0} is $$1.24533{}92504{}91216{}63010{}69081{}87858{}37457{}65950,$$
for which Mathematica reports an accuracy of 37~digits.

While~\eqref{eqn:beta0} makes it appear that Magma's estimate of 4000 terms is conservative,
this is partially explained by the fact that we are evaluating
at the point $s=1$.  If we wanted to calculate $L(\frac12 + 100 i)$
then many more terms would be needed.  Nevertheless, we see that
180 terms, or even 1000 terms, are not adequate for the high precision
evaluation we seek.

As we will explain in Section~\ref{sec:averaging},
we can achieve high precision by evaluating the L-function several times,
and then taking a linear combination of those evaluations.
Thus, we need to evaluate~\eqref{eqn:genbeta} for other values of $\beta$,
say
$\beta=\frac32$:
\begin{align}\label{eqn:beta3half}
L_{\frac32}(1) =\mathstrut&  1.870\, b_1 + 0.937\, b_2 + 0.017\, b_3 + \cdots +
0.0097\, b_{17}  \notag \\
&+ \cdots -
2.10\times 10^{-8}\, b_{101}  +
\cdots -
9.44\times 10^{-12}\, b_{181} \notag \\
&+
\cdots
+
4.6\times 10^{-19}\, b_{499}
+
\cdots
+
4.7\times 10^{-25}\, b_{1009} \notag \\
& +
\cdots
+
1.3 \times 10^{-29}\, b_{1499}
+
\cdots
-
4.2 \times 10^{-33}\, b_{1999} \notag \\
&
+
\cdots + 5.2 \times 10^{-47}\, b_{4999} + \cdots.
\end{align}

The choice of $\beta=\frac32$ seems worse, because the contributions of
the coefficients are decreasing less rapidly, so more terms will
be required in order to obtain a given precision.  This is indeed true,
for among test functions of this form,
$\beta=0$ has the contributions decreasing most rapidly.

From \eqref{eqn:beta0} or \eqref{eqn:beta3half} we can determine a value for $L(1)$
by using the known coefficients
and estimating the others with the Ramanujan bound $|b_p|\le 8$.
Note that the known coefficients include not only $b_n$ for $n\le 180$,
but also some larger numbers such as $b_{875} = b_7 b_{125}$.
For $\beta=0$ we find
\begin{align}\label{eqn:beta0est}
L_0(1) =\mathstrut& 1.798902826118503606167865 + 8.56 \times 10^{-14}\, b_{181}   \notag \\
 &+ 3.81\times 10^{-14}\, b_{191} + 3.25 \times 10^{-14}\, b_{193}
+ 2.37 \times 10^{-14}\, b_{197}\notag \\
&+\cdots + 1.1 \times 10^{-21}\, b_{499}
 +
\cdots
+
5.5\times 10^{-29}\, b_{1009}\notag \\
&
+
\cdots
+
7.3 \times 10^{-34}\, b_{1499}
+
\cdots
+
8.3 \times 10^{-38}\, b_{1999}   \notag \\
&+
\cdots + 6.7 \times 10^{-53} \, b_{4999} + \cdots  \notag \\
=\mathstrut& 1.79890 28261 18503  \pm 1.8 \times 10^{-12},
\end{align}
and for $\beta=\frac32$:
\begin{align}\label{eqn:beta3halfest}
L_{\frac32}(1) =\mathstrut& 1.798902826123555372082651 + 9.44 \times 10^{-12}\, b_{181}  \notag \\
 &+ 9.42 \times 10^{-12}\, b_{191}  + 8.85 \times 10^{-12}\, b_{193}
+ 7.54 \times 10^{-12}\, b_{197}\notag \\
 &+\cdots + 4.6 \times 10^{-19}\, b_{499}
+
\cdots
+
4.7\times 10^{-25}\, b_{1009} \notag \\
&
+
\cdots
+
1.3 \times 10^{-29}\, b_{1499}
+
\cdots+
4.2 \times 10^{-33}\, b_{1999}  \notag \\
&+
\cdots + 5.2 \times 10^{-47} b_{4999} + \cdots   \notag \\
=\mathstrut& 1.79890 28261 235 \pm 5.5\times 10^{-10} .
\end{align}
Those values are far from our goal of 30 decimal digits of accuracy.

In \eqref{eqn:beta0est} and \eqref{eqn:beta3halfest}, and below,
the expression
$a=b \pm c$ means that the true value of $a$ lies in the interval
$(b-c, b+c)$.
To estimate the error, we used the first $20,000$ Dirichlet coefficients.

\subsection{Averaging is better}\label{sec:averaging}

In this section we use the main idea of~\cite{FR}, which is that
one can obtain
a more precise evaluation of $L(1)$ by combining
$L_\beta(1)$ for several auxiliary functions.  That is, if
\begin{equation}\label{eqn:sum1}
\sum w_j = 1,
\end{equation}
then
\begin{equation}\label{eqn:generalavg}
L(1) = \sum  w_j L_{\beta_j}(1) .
\end{equation}
As we will see, suitable choices of $w_j$ and $\beta_j$
in \eqref{eqn:generalavg} will give a value of $L(1)$
with a small error term.
For example, consider
\begin{align}\label{eqn:fourweights}
 \{\beta_j\} = \mathstrut& \{0,\  \tfrac12,\  1,\  \tfrac32 \}  \\
\{w_j\} = \mathstrut& \{5.595844269,\  -5.074113323,\  0.484231975,\  -0.0059629212\}. \notag
\end{align}
The ``magic numbers'' in \eqref{eqn:fourweights} were chosen (using a least-squares fit)
to minimize the contribution of $b_{181},  b_{191}, \ldots$ in \eqref{eqn:avg4} below.

Using those values for $w_j$ and $\beta_j$, we find:
\begin{align}\label{eqn:avg4}
L(1)  =\mathstrut& \sum  w_j L_{\beta_j}(1) \notag \\
=\mathstrut& 1.798902826118603393418629 + 1.34 \times 10^{-17}\, b_{181}   \notag \\
 &+ 4.96 \times 10^{-18}\, b_{191}+ 3.49 \times 10^{-17}\, b_{193}
+ 5.01 \times 10^{-17}\, b_{197} \notag \\
&+\cdots + 1.1 \times 10^{-20}\, b_{499} +
\cdots
+
1.9 \times 10^{-27}\, b_{1009}  \notag \\
&
+ \cdots+9.4 \times 10^{-32}\, b_{1499} + \cdots +
3.1 \times 10^{-35}\, b_{1999}\notag \\
&
+ \cdots -3.1 \times 10^{-49} \, b_{4999} +\cdots \notag \\
= \mathstrut& 1.79890 28261 18603 39 \pm 3 \times 10^{-14} .
\end{align}
The error in \eqref{eqn:avg4} compared to \eqref{eqn:beta0est} and \eqref{eqn:beta3halfest}
should be somewhat surprising.
The error has decreased by a
factor of 60 by combining only 4 evaluations of~$L(1)$.
Suppose that, instead, we wanted to
improve the error in \eqref{eqn:beta0est} by determining more coefficients.
We would need to determine the value of $b_p$ for the 11 primes $181 \le p \le 239$
in order to have a comparable decrease in the error term.

The improved error in \eqref{eqn:avg4} indicates that the contributions from the
unknown coefficients, for different auxiliary functions,
are negatively correlated.
The exact nature of this correlation has not been described analytically:
we take it as an empirical fact.

By combining the evaluation of $L_\beta(1)$ for
$\beta \in \{0, 1/10, 2/10, \ldots, 30/10\}$, with suitably
chosen (by least-squares) ``weights'' $\{w_0, w_1, \ldots, w_{30}\}$, we find
\begin{equation}
L(1) = 1.79890 28261 18603 03245 57227 72619 \pm 6 \times 10^{-26}
\end{equation}
and in the same way,
\begin{equation}
L(3) = 1.10545 68879 51321 63036 93593 41690  \pm 3 \times 10^{-27} .
\end{equation}
Thus
\begin{align}\label{eqn:manyzeros}
\frac{\pi^8 L(1)}{L(3)} =\mathstrut &
       15440.62500 00000 00000 00000 000096
\pm 6 \times 10^{-22} \cr
= \mathstrut & \frac{3^4 \cdot 5^2 \cdot 61}{2^3} \pm 6 \times 10^{-22} .
\end{align}
This precision of 27 digits does not quite meet our goal of 30 digits of accuracy in the
final result, but the identification of ${\pi^8 L(1)}/{L(3)}$ as a rational
number seems convincing.  Similar calculations produced the other
values in Section~\ref{compres}.

Note that this approach to evaluating L-functions requires significantly
more computation than the methods used when more coefficients are known.
With known coefficients, the sum in \eqref{eqn:appfe} and the integrals
in \eqref{eqn:mellin} can be interchanged, so only two integrals need to be
computed numerically.  For the above calculation,
it was necessary to compute each of the thousands of integrals separately.
Furthermore, those integrals were computed multiple times:
once for each different auxiliary function.  Only after evaluating all those
integrals could we optimally combine them to minimize the
contributions of the unknown coefficients.   If the appropriate combinations
could somehow be determined in advance, our ability to evaluate higher-degree L-functions
would be substantially improved.

\section{Expected consequences of congruences revisited: the Bloch-Kato conjecture}
\subsection{Statement of the conjecture}
Recall the rank-$8$ motive $M=M_f\otimes M_F$ such that $L(M,s)=L(f\otimes F,s)$. (We shall assume at least the existence of a premotivic structure comprising realisations and comparison isomorphisms, as defined in \cite[1.1.1]{DFG}.) In our examples the coefficient field is $\QQ$. Let $q>j+2k+\ell-3$ be a prime number. Choose a $\ZZ_{(q)}$-lattice $T_B$ in the Betti realisation $M_B$ in such a way that $T_{q}:=T_B\otimes \ZZ_q$ is a $\Gal(\Qbar/\Q)$-invariant lattice in the $q$-adic realisation $M_q$, via the comparison isomorphism $M_B\otimes\Q_q\simeq M_q$. Then choose a $\ZZ_{(q)}$-lattice $T_{\dR}$ in the de Rham realisation
 $M_{\dR}$ in such a way that
$$\VVV(T_{\dR}\otimes \ZZ_q)=T_{q}$$ as $\Gal(\Qbar_q/\QQ_q)$-representations, where $\VVV$ is the version of the Fontaine-Lafaille functor used in \cite{DFG}. Since $\VVV$ only applies to filtered $\phi$-modules, where $\phi$ is the crystalline Frobenius, $T_{\dR}\otimes \ZZ_q$ must be $\phi$-stable. Anyway, this choice ensures that the $q$-part of the Tamagawa factor at $q$ is trivial (by \cite[Theorem 4.1(iii)]{BK}), thus simplifying the Bloch-Kato conjecture below. The condition $q>j+2k+\ell-3$ ensures that the condition (*) in \cite[Theorem 4.1(iii)]{BK} holds.

Let $t$ be a critical point at which we evaluate the $L$-function.
Let $M(t)$ be the corresponding Tate twist of the motive.
Let $\Omega(t)$ be a Deligne period scaled according to the above choice, i.e. the determinant of the isomorphism
$$M(t)_B^+\otimes\CC\simeq (M(t)_{\dR}/\Fil^0)\otimes\CC,$$
calculated with respect to bases of $(2\pi i)^tT_B^{(-1)^t}$ and $T_{\dR}/\Fil^t$, so well-defined up to $\ZZ_{(q)}^{\times}$.

The following formulation of the ($q$-part of the) Bloch-Kato conjecture, as applied to this situation, is based on \cite[(59)]{DFG} (where $\Sigma$ was non-empty, though), using the exact sequence in their Lemma 2.1.
\begin{conj}[Bloch-Kato]\label{BK}
$$\ord_{q}\left(\frac{L(M,t)}{\Omega(t)}\right)$$
$$=\ord_{q}\left(\frac{\#H^1_{f}(\Q,T_{q}^*(1-t)\otimes (\QQ_q/\ZZ_q))}{\#H^0(\Q,T_{q}^*(1-t)\otimes(\QQ_q/\ZZ_q))\#H^0(\Q,T_{q}(t)\otimes(\QQ_q/\ZZ_q))}\right).$$
\end{conj}
Here, $T_{q}^*=\Hom_{\ZZ_q}(T_{q},\ZZ_q)$, with the dual action of $\Gal(\Qbar/\Q)$. This is an invariant $\ZZ_q$-lattice in $M_{q}^*\simeq M_q(j+2k+\ell-4)$, so $T_q^*(1-t)$ is a lattice in $M_q(j+2k+\ell-3-t)$. On the right hand side, in the numerator, is a Bloch-Kato Selmer group with local conditions (unramified at $p\neq q$, crystalline at $p=q$) for all finite primes $p$.
\subsection{Global torsion and Kurokawa-Mizumoto type congruences}
We revisit the situation of \S 3.1. Recall that $\lambda_F(p)$ denotes the eigenvalue of the genus-$2$ Hecke operator $T(p)$ acting on the cuspidal eigenform $F$. The $q$-adic realisations $M_{f,q}$ and $M_{F,q}$ should be $2$-and $4$-dimensional $\QQ_q$ vector spaces with continuous linear actions $\rho_f$, $\rho_F$ of $\Gal(\Qbar/\QQ)$, crystalline at $q$, unramified at all primes $p\neq q$. For primes $p\neq q$, we should have $$a_p(f)=\Tr(\rho_f(\Frob_p^{-1}))\,\,\,\,\,\text{and}\,\,\,\,\,\lambda_F(p)=\Tr(\rho_F(\Frob_p^{-1})).$$ Galois representations with these properties are known to exist, by theorems of Deligne and Weissauer \cite{De2,We}. By Poincar\'e duality, $M_{f,q}^*\simeq M_{f,q}(\ell-1)$ and $M_{F,q}^*\simeq M_{F,q}(j+2k-3)$. Choosing $\Gal(\Qbar/\QQ)$-invariant $\ZZ_q$-lattices in $M_{f,q}$ and $M_{F,q}$, then reducing mod $q$, we obtain residual representations $\rhobar_f$ and $\rhobar_F$. We suppose that (as in Example 1) $\rhobar_f$ is irreducible, in which case it is independent of the choice of lattice in $M_{f,q}$. The congruence
$$\lambda_F(p)\equiv a_p(f)(1+p^{k-2})\pmod{q},$$
interpreted as a congruence of traces of Frobenius, implies that the composition factors of $\rhobar_F$ are $\rhobar_f$ and $\rhobar_f(2-k)$. Which is a submodule and which is a quotient will depend on the choice of lattice in $M_{F,q}$.

Looking at the denominator of the Bloch-Kato formula, with $T_q$ the tensor product of the $\ZZ_q$-lattices referred to above, on which $\Gal(\Qbar/\QQ)$ acts by $\rho_f\otimes\rho_F\simeq \rho_f^*(1-\ell)\otimes\rho_F$, the $q$-torsion in $H^0(\Q,T_{q}(t)\otimes(\QQ_q/\ZZ_q))$ is $(\rhobar_f^*\otimes\rhobar_F(t+1-\ell))^{\Gal(\Qbar/\QQ)}$, which is $\Hom_{\FF_q[\Gal(\Qbar/\QQ)]}(\rhobar_f,\rhobar_F(t+1-\ell))$. Recalling that $\ell=j+k$, this is the same as $\Hom_{\FF_q[\Gal(\Qbar/\QQ)]}(\rhobar_f(2-k),\rhobar_F(t+3-j-2k))$. This can be non-trivial only for $t\equiv \ell-1\pmod{q-1}$ (if $\rho_f$ is a submodule of $\rho_F$) or for $t\equiv j+2k-3\pmod{q-1}$ (if $\rho_f(2-k)$ is a submodule of $\rho_F$). The only such $t$ in the critical range $\ell\leq t\leq j+2k-3$ (using $q>j+2k+\ell-4$ from \S 5.1, or even just $q>2\ell$ from \S 3.1) is $t=j+2k-3$. So, with a suitable choice of lattice, and $t=j+2k-3$, we can have a factor of $q$ in the denominator of the conjectural formula for $\frac{L(M,t)}{\Omega(t)}$, which appears to provide some explanation for the $q$ in the denominator of $\frac{L(j+2k-3,f\otimes F)}{\pi^8L(j+2k-5,f\otimes F)},$ observed in Example 1.

Note that the factor $L(s-(k-2),\Sym^2 f)$ in the expression in \S 3.1 for $L(f\otimes[f]_j)$ has trivial zeros at the points $s=\ell, \ell+2$ paired with $s=j+2k-3, j+2k-5$ by the functional equation. This is because $\ell-(k-2)=j+2$ and $\ell+2-(k-2)=j+4$, which are even and in the range $1\leq t\leq \ell-1$, at least if $k>4$. This suggests that the orders of Selmer groups may contribute cancelling factors of $q$ to the numerators of $L_{\alg}(j+2k-3,f\otimes F)$ and $L_{\alg}(j+2k-5,f\otimes F)$, something we overlooked in the final paragraph of \cite[\S 4.2]{DHR}. Note also that one can make a similar construction of global torsion elements with respect to congruences of ``Yoshida type'' (which appear in \cite[Conjecture 10.7]{BFvdG}), but in that case there are no critical $L$-values.
\subsection{Moving between Selmer groups via Harder type congruences}
Now we revisit the situation of \S 3.2. There $t=(\ell/2)+j+k-1$, and $T_q^*(1-t)$ is a lattice in $M_q(j+2k+\ell-3-t)=M_q((\ell/2)+k-2)$.
By an analogue of the Birch and Swinnerton-Dyer
conjecture, vanishing of $L(f,\ell/2)$ should suffice for the
non-triviality of $H^1_f(\QQ, M_{f,q}(\ell/2))$ (again defined using local conditions). (See the
``conjectures'' $C_r(M)$  in \S 1 of \cite{Fo}, and
$C^i_{\lambda}(M)$ in \S 6.5 of \cite{Fo}.) The sign in the functional equation of $L(f,s)$ is $(-1)^{\ell/2}=-1$,
so the parity of the order of vanishing at $s=\ell/2$ is odd. Assuming that $\rhobar_f$ is irreducible, the
conditions of \cite[Theorem B]{N2} are satisfied. Hence $H^1_f(\QQ,
M_{f,q}(\ell/2))$ is non-trivial (because the parity of its rank
is also odd). If we were to impose a
condition that $f$ is ordinary at $q$ (i.e. $q\nmid
a_{q}(f)$), then we could alternatively get this from either \cite[Th\'eor\`eme A]{SU} or the main theorem of \cite[\S 12]{N1}.

Anyway, from this one easily obtains a non-zero element $$c''\in H^1(\QQ,\rhobar_f(\ell/2)).$$ Assuming irreducibility of $\rhobar_f$, it is a consequence of Harder's conjectured congruence that the composition factors of $\rhobar_F$ are $\rhobar_g$, $\FF_q(2-k)$ and $\FF_q(1-j-k)$. If we choose the $\Gal(\Qbar/\QQ)$-invariant $\ZZ_q$-lattice in $M_{F,q}$ in such a way that the composition factor $\FF_q(2-k)$ of $\rhobar_F$ is a submodule, then $\rhobar_f(2-k)$ is a submodule of $\rhobar_f\otimes\rhobar_F$, so $\rhobar_f(\ell/2)$ is a submodule of $\rhobar_f\otimes\rhobar_F((\ell/2)+k-2)$. Thus we may map $c''$ to $H^1(\QQ,\rhobar_f\otimes\rhobar_F((\ell/2)+k-2))$, thence to $H^1(\Q,T_{q}^*(1-t)\otimes (\QQ_q/\ZZ_q))$. Assuming that $\rhobar_f\not\simeq\rhobar_g$ (e.g. if $\ell\neq k'$ and $q>\max\{\ell,k'\}$) one easily checks that $H^0(\QQ,\rhobar_f\otimes\rhobar_F((\ell/2)+k-2))$ is trivial, from which it follows that this produces a {\em non-zero} element of $H^1(\Q,T_{q}^*(1-t)\otimes (\QQ_q/\ZZ_q))$. If $q>j+2k+\ell-3$ one can in fact show, as in the proof of \cite[Proposition 5.1]{DIK}, that this element is in $H^1_f(\Q,T_{q}^*(1-t)\otimes (\QQ_q/\ZZ_q))$. This puts a factor of $q$ in the numerator of the conjectural formula for $\frac{L(M,t)}{\Omega(t)}$, which appears to provide some explanation for the $q$ in the numerator of $\frac{L((\ell/2)+j+k-1,f\otimes F)}{\pi^8L((\ell/2)+j+k-3,f\otimes F)}$, observed in Example 2.

Analogous situations were already considered in \cite[\S 8,\S 11,\S 14]{Du3} and \cite[Conjecture 5.4, Corollary 8.6]{DIK}.

\section{Endoscopic congruences for $\SO(7)$.}\label{endo7}
{\bf Example 3. }
Recall that when $\ell=16$ and $(j,k)=(4,12)$, we found an apparent factor of $71$ in the numerator of
$\frac{\pi^8L(1,\pi_f\otimes \pi_F)}{L(3,\pi_f\otimes \pi_F)}=\frac{\pi^8L((\ell+j+2k-2)/2,f\otimes F)}{L((\ell+j+2k+2)/2,f\otimes F)}.$
With $q=71$ and $t=(\ell+j+2k-2)/2$ (which is the integer immediately to the right of the centre of the functional equation), we would like to construct a non-zero element in $H^1_f(\Q,T_{q}^*(1-t)\otimes (\QQ_q/\ZZ_q))$ to explain this. The $q$-torsion in $T_{q}^*(1-t)\otimes (\QQ_q/\ZZ_q)$ is (the space of) $\rhobar_f\otimes\rhobar_F((\ell+j+2k-4)/2)\simeq \Hom_{\FF_q}(\rhobar_f((\ell+2-j-2k)/2),\rhobar_F)$, using $\rhobar_f^*\simeq \rhobar_f(\ell-1)$. Note that the Hodge type of $M_f((\ell+2-j-2k)/2)$ is $\{((j+2k-\ell-2)/2,\ell-1+(j+2k-\ell-2)/2),(\ell-1+(j+2k-\ell-2)/2,(j+2k-\ell-2)/2)\}$, which is $\{(5,20),(20,5)\}$ in this case, and the effect of the twist is to make the ``weight'' $w=$``$p+q$''$=\ell-1+(j+2k-\ell-2)=j+2k-3$, equal to that of the Hodge type $\{(0,j+2k-3),(k-2,j+k-1),(j+k-1,k-2),(j+2k-3,0)\}$ of $M_F$, which is $\{(0,25),(10,15),(15,10),(25,0)\}$ in this case.
This raises the possibility that $\rhobar_f((\ell+2-j-2k)/2)=\rhobar_f(-5)$ and $\rhobar_F$ could both occur as composition factors in the reduction mod $q$ of an invariant $\ZZ_q$-lattice in a $6$-dimensional Galois representation coming from the $q$-adic realisation of a rank $6$ motive, pure of weight $j+2k-3$, with Hodge type the union of those of $M_f((\ell+2-j-2k)/2)$ and $M_F$.

If such a $6$-dimensional $q$-adic Galois representation $\tilde{\rho}$ exists, and if it is irreducible, then adapting a well-known construction of Ribet \cite{Ri}, there exists an invariant $\ZZ_q$-lattice whose reduction mod $q$ is a non-trivial extension of $\rhobar_f((\ell+2-j-2k)/2)$ by $\rhobar_F$ (both of which we suppose to be irreducible). This gives a non-zero element in $H^1(\QQ,\Hom_{\FF_q}(\rhobar_f((\ell+2-j-2k)/2),\rhobar_F))$. Since $H^0(\QQ,\Hom_{\FF_q}(\rhobar_f((\ell+2-j-2k)/2),\rhobar_F))$ is trivial, the image of this element in $H^1(\Q,T_{q}^*(1-t)\otimes (\QQ_q/\ZZ_q))$ is non-zero. If we also suppose that $\tilde{\rho}$ is unramified at all $p\neq q$, crystalline at $q$, then one can show (using $q>(3j+6k+\ell-8)/2$) that it lies in fact in $H^1_f(\Q,T_{q}^*(1-t)\otimes (\QQ_q/\ZZ_q))$, as desired.
It remains to explain where $\tilde{\rho}$ comes from.

Given a self-dual, cuspidal, automorphic representation $\pi$ of $\GL_6(\A)$, there is an associated $\tilde{\rho}:\Gal(\Qbar/\QQ)\rightarrow \GL_6(\QQ_q)$ (see \cite[Remark 7.6]{Sh}). If $\pi$ is unramified at all finite places then $\tilde{\rho}$ is unramified at all $p\neq q$ and crystalline at $q$. It is not currently known to be irreducible, but we shall assume that, as expected, it is, so that the above construction applies. The infinitesimal character of $\pi_{\infty}$ determines the Hodge-Type of the conjectural motive of which $\tilde{\rho}$ should be the $q$-adic realisation \cite{Cl} (and the Hodge-Tate weights of $\tilde{\rho}|_{\Gal(\Qbar_q/\QQ_q)}$). A self-dual, cuspidal automorphic representation of $\PGL_6(\A)$ discovered by Chenevier and Renard, denoted $\Delta_{25,15,5}$ in \cite[Table 13]{CR}, has the correct infinitesimal character. By Arthur's symplectic-orthogonal alternative \cite[Theorem* 3.9]{CR}, it is the functorial lift of a discrete automorphic representation of $\SO(4,3)(\A)$.

Let $\SO(7)$ be the special orthogonal group of the $E_7$ root lattice, the even, positive-definite lattice of discriminant $2$, unique up to isomorphism.
This is a semi-simple group over $\ZZ$, and $\SO(7)(\ZZ)\simeq W(E_7)^+$, the even subgroup of the Weyl group, of order $1451520$. For $\mu=a_1e_1+a_2e_2+a_3e_3$ (in the notation of \cite[5.2]{CR}), with $a_1,a_2,a_3\in\ZZ$ and $a_1\geq a_2\geq a_3\geq 0$, let $V_{\mu}$ be the space of the complex representation $\theta_{\mu}$ of $\SO(7)$ with highest weight $\mu$, and let $\rho:=\frac{5}{2}e_1+\frac{3}{2}e_2+\frac{1}{2}e_3$. The infinitesimal character of the representation $\theta_{\mu}$ of $\SO(7)(\RR)$ is $\mu+\rho$. Let $K$ be the open compact subgroup $\prod_p \SO(7)(\ZZ_p)$ of $\SO(7)(\A_f)$, and let
\begin{multline*}M(V_{\mu},K):=\{f:\SO(7)(\A_f)\rightarrow V_{\mu}: f(gk)=f(g)~\forall k\in K, \\f(\gamma g)=\theta_{\mu}(\gamma)(f(g))~\forall \gamma\in \SO(7)(\QQ)\}\end{multline*}
be the space of $V_{\mu}$-valued algebraic modular forms with level $K$ (i.e. ``level $1$''), where $\A_f$ is the ``finite'' part of the adele ring.
Since $\#(\SO(7)(\QQ)\backslash \SO(7)(\A_f)/K)=1,$ $M(V_{\mu},K)$ can be identified with the fixed subspace $V_{\mu}^{\SO(7)(\ZZ)}$.

If we let $\mu=10e_1+6e_2+2e_3$, so that $\mu+\rho=(1/2)(25e_1+15e_2+5e_3)$, then $M(V_{\mu},K)$ is $2$-dimensional, spanned by $K$-fixed vectors of automorphic representations of $\SO(7)(\A)$ whose Arthur parameters are $\Delta_{25,15,5}$ and an endoscopic lift denoted $\Delta_{25,5}\oplus\Delta_{15}$ in \cite[Table 13]{CR}. Note that $(25,5)=(j+2k-3,j+1)$ and $15=\ell-1$, with $(j,k)=(4,12)$ and $\ell=16$. M\'egarban\'e has calculated the traces of certain Hecke operators $T(p)$ (for $p\leq 53$) on spaces including this one \cite{Me3}. The contribution
$p^5a_p(f)+\lambda_F(p)$ from $\Delta_{25,5}\oplus\Delta_{15}$ is easily subtracted off (as below) to find the eigenvalue denoted $T(p)(\Delta_{25,15,5})$, in fact M\'egarban\'e has recorded the results in \cite[Tables 2,4,5,6]{Me}. A congruence of Hecke eigenvalues
$$T(p)(\Delta_{25,15,5})\equiv p^5a_p(f)+\lambda_F(p)\pmod{71}$$
for all primes $p$ would, viewing them as traces of Frobenius, imply that $\overline{\tilde{\rho}}$ has composition factors $\rhobar_f((\ell+2-j-2k)/2)$ and $\rhobar_F$, which is what we need. We confirmed this congruence for all $p\leq 53$. In the tables below, we show the results for $2\leq p\leq 11$ and $p=53$.
\vskip10pt
\begin{tabular}{|c|c|c|c|}\hline $p$ & $a_p(f)$ & $\lambda_F(p)$ & $\Tr(T(p)|V_{\mu}^{\SO(7)(\ZZ)})$\\\hline $2$ & $216$ & $-96$ & $6816$\\$3$ & $-3348$ & $-527688$ & $-474120$\\$5$ & $52110$ & $596139180$ & $145932324$\\$7$ & $2822456$ & $-3608884496$ & $49205357040$\\$11$ & $20586852$ & $3047542095144$ & $3229012641000$\\$53$ & $6797151655902$ & $-3921035060705523617268$ & $-8934610079$\\ & & & $5036491708$\\\hline
\end{tabular}
\vskip10pt
\begin{tabular}{|c|c|c|}\hline $p$ & $T(p)(\Delta_{25,15,5})$ & $-T(p)(\Delta_{25,15,5})+p^5a_p(f)+\lambda_F(p)$\\\hline $2$ & $0$ & $2^5.3.\mathbf{71}$\\$3$ & $867132$ & $-2^7.3^5.\mathbf{71}$\\$5$ & $-613050606$ & $2^9.3.23.\mathbf{71}.547$\\$7$ & $5377223544$ & $2^8.3^4.7^2.13.41.\mathbf{71}$\\$11$ & $-3134062555596$ & $2^7.3.23.\mathbf{71}.15145211$\\$53$ & $989150772174783875874$ & $-2^9.3^3.\mathbf{71}.73.1031.27990002153$\\\hline
\end{tabular}
\vskip10pt
Recall from \S \ref{compres} that when $\ell=18$ and $(j,k)=(4,15)$, we see a factor of $193$ in the numerator of
$\frac{\pi^8L(1,\pi_f\otimes \pi_F)}{L(3,\pi_f\otimes \pi_F)}.$
Were we to try to account for this in a similar manner, we would have to look at $\mu+\rho=(1/2)(31e_1+17e_2+5e_3)$, in which case $M(V_{\mu},K)$ is $8$-dimensional and $\Delta_{31,17,5}$ is not unique (using \cite{CRtab}). Hence we cannot extract Hecke eigenvalues from traces of Hecke operators in the same way, so are unable to test the expected mod $193$ congruence.

{\bf Example 4. }
Recall that when $\ell=16$ and $(j,k)=(6,10)$, we found an apparent factor of $61$ in the numerator of
$\frac{\pi^8L(1,\pi_f\otimes \pi_F)}{L(3,\pi_f\otimes \pi_F)}=\frac{\pi^8L((\ell+j+2k-2)/2,f\otimes F)}{L((\ell+j+2k+2)/2,f\otimes F)}.$
In the same way as for the previous example, this could be explained by a certain mod $61$ congruence of Hecke eigenvalues of algebraic modular forms for $\SO(7)$, this time with $\mu=9e_1+6e_2+3e_3$. We have again used M\'egarban\'e's data to verify the congruence for $p\leq 53$.
Similarly, it is necessary to subtract from the trace a contribution from $\Delta_{23,7}\oplus\Delta_{15}$, to get the Hecke eigenvalue for $\Delta_{23,15,7}$.
\vskip10pt
\begin{tabular}{|c|c|c|}\hline $p$ &  $\lambda_F(p)$ & $\Tr(T(p)|V_{\mu}^{\SO(7)(\ZZ)})$\\\hline $2$ & $1680$ & $4416$\\$3$ & $-6120$ & $148104$\\$5$ & $2718300$ & $-89271276$\\$7$ & $6916898800$ & $10652657232$\\$11$ & $-1417797110136$ & $-764339838888$\\$53$ & $-15111411349636553220$ & $86535126376033794804$\\\hline
\end{tabular}
\vskip10pt
\begin{tabular}{|c|c|c|}\hline $p$ & $T(p)(\Delta_{23,15,7})$ & $-T(p)(\Delta_{23,15,7})+p^4a_p(f)+\lambda_F(p)$\\\hline $2$ & $-720$ & $2^5.3.\mathbf{61}$\\$3$ & $425412$ & $-2^8.3^2.5.\mathbf{61}$\\$5$ & $-124558326$ & $2^{10}.3.\mathbf{61}.853$\\$7$ & $-3040958424$ & $2^9.3^7.5.7^2.\mathbf{61}$\\$11$ & $352045171116$ & $-2^8.3.5.7^2.31.\mathbf{61}.4127$\\$53$ & $48013741730657079162$ & $-2^{10}.3^3.5.17.\mathbf{6
1}.66215793179$\\\hline
\end{tabular}
\vskip10pt
In fact, M\'egarban\'e has very recently proved this congruence unconditionally for all $p$ \cite[Th\'eor\`eme 1.0.3(x)]{Me1}.
He uses scalar-valued algebraic modular forms for $\mathrm{SO}(25)$, in the manner of Chenevier and Lannes's proof of Harder's mod $41$ congruence using $\mathrm{O}(24)$ \cite[Chapter 10, Theorem* 4.4(1)]{CL}. He found that the modulus of the congruence is in fact $5856=2^5.3.61$.
\section{An Eisenstein congruence for $\SO(9)$.}\label{eis9}
{\bf Example 5. }
Recall that when $\ell=16$ and $(j,k)=(4,12)$, we found an apparent factor of $17$ in the denominator of
$\frac{\pi^8L(1,\pi_f\otimes \pi_F)}{L(3,\pi_f\otimes \pi_F)}=\frac{\pi^8L((\ell+j+2k-2)/2,f\otimes F)}{L((\ell+j+2k+2)/2,f\otimes F)},$
and in the numerator of $\frac{\pi^8L(3,\pi_f\otimes \pi_F)}{L(5,\pi_f\otimes \pi_F)}$, so apparently in the numerator of $L_{\alg}((\ell+j+2k+2)/2,f\otimes F)$. With $q=17$ and $t=(\ell+j+2k+2)/2$ (which is no longer immediately to the right of the centre of the functional equation), we would like to construct a non-zero element in $H^1_f(\Q,T_{q}^*(1-t)\otimes (\QQ_q/\ZZ_q))$ to try to explain this, though the condition $q>j+2k+\ell-3$ does not hold here. The $q$-torsion in $T_{q}^*(1-t)\otimes (\QQ_q/\ZZ_q)$ is (the space of) $\rhobar_f\otimes\rhobar_F((\ell+j+2k-8)/2)\simeq \Hom_{\FF_q}(\rhobar_f((\ell+6-j-2k)/2),\rhobar_F)$, using $\rhobar_f^*\simeq \rhobar_f(\ell-1)$.

We would like to see $\rhobar_f((\ell+6-j-2k)/2)=\rhobar_f(-3)$ and $\rhobar_F$ both occurring as composition factors in the reduction mod $q$ of an invariant $\ZZ_q$-lattice in a $6$-dimensional Galois representation coming from the $q$-adic realisation of a rank $6$ motive. Then we could apply the construction of Ribet again (though $q$ is not large enough now for us to prove the local condition at $q$). The problem is, $M_F$ still has Hodge type $\{(0,25),(10,15),(15,10),(25,0)\}$, of weight $25$, while the Hodge type of $M_f(-3)$ is $\{(3,18),(18,3)\}$, of weight only $21$.
What we need to do is to balance $\rhobar_f(-3)$ with another composition factor $\rhobar_f(-7)$, noting that the Hodge type of $M_f(-7)$ is $\{(7,22),(22,7)\}$, and $3+22=7+18=25$. Now $\rhobar_f(-3)$, $\rhobar_f(-7)$ and $\rhobar_F$ could all occur as composition factors in the reduction mod $q$ of an invariant $\ZZ_q$-lattice in an $8$-dimensional Galois representation coming from the $q$-adic realisation of a rank $8$ motive, pure of weight $25$, with Hodge type $\{(0,25),(3,22),(7,18),(10,15),(15,10),(18,7),(22,3),(25,0)\}$. Although this is not the union of the Hodge types of $M_f(-3), M_f(-7)$ and $M_F$, the union of the sets of Hodge-Tate weights of their $q$-adic realisations restricted to $\Gal(\Qbar_q/\QQ_q)$ is  $\{0,3,7,10,15,18,22,25\}$.

This time a self-dual, cuspidal, automorphic representation of $\PGL_8(\A)$ discovered by Chenevier and Renard  \cite[Corollary**6.5, Table 8]{CR}, denoted $\Delta_{25,19,11,5}$ in their notation, has the correct infinitesimal character. By Arthur's symplectic-orthogonal alternative \cite[Theorem* 3.9]{CR}, it is the functorial lift of a discrete automorphic representation of $\SO(5,4)(\A)$. Again, there is an associated $\tilde{\rho}:\Gal(\Qbar/\QQ)\rightarrow \GL_8(\QQ_q)$ (see \cite[Remark 7.6]{Sh}). The relevant space of algebraic modular forms for $\SO(9)$ is $3$-dimensional, spanned by Hecke eigenforms that are vectors in automorphic representations of $\SO(9)(\A)$ with Arthur parameters $\Delta_{25,19,11,5}$ and $\Delta^2_{25,19,5}\oplus\Delta_{11}$, \cite[Table 1]{Me}. Here $\Delta^2_{25,19,5}$ stands for a pair of self-dual, cuspidal, automorphic representations of $\PGL_6(\A)$, and $\Delta^2_{25,19,5}\oplus\Delta_{11}$ for a pair of endoscopic lifts. To get the Hecke eigenvalues we want, for $\Delta_{25,19,11,5}$, one must subtract the endoscopic contributions from the traces computed by M\'egarban\'e for $p\leq 7$ \cite{Me3}. Also, computing the trace of the $\SO(7)$ $T(p)$ on $\Delta^2_{25,19,5}$ similarly requires the subtraction of an endoscopic contribution by $\Delta_{25,5}\oplus\Delta_{19}$ from a trace on a whole space of algebraic modular forms. One can obtain the $T(p)(\Delta_{25,19,11,5})$ directly from \cite[Table 7]{Me}.

The congruence verified in the second table for $p\leq 7$, if it held for all $p$, would imply that (with $q=17$) $\tilde{\rhobar}$ has composition factors $\rhobar_f(-3)$, $\rhobar_f(-7)$ and $\rhobar_F$. Note that $\rhobar_f$ is certainly irreducible, by \cite[Corollary to Theorem 4]{SwD}, and the irreducibility of $\rhobar_F$ can presumably be checked as in \cite[Proposition 4.10]{CL}. This congruence (disregarding the smallness of $q$) is an instance of the kind considered in \cite{BD}, in the case $G=\SO(5,4)$. The expression $(p^3+p^7)a_p(f)+\lambda_F(p)$ is the eigenvalue of $T(p)$ on an automorphic representation of $G(\A)$ induced from a maximal parabolic subgroup with Levi subgroup $M\simeq\GL(2)\times\SO(3,2)$. (Harder's congruences above are a different instance, as explained in \cite[\S 7]{BD} and \cite[\S 3.2]{BDM}.)
\vskip10pt
\begin{tabular}{|c|c|c|c|c|}\hline $p$ & $\Tr(T(p)|V_{\mu}^{\SO(7)(\ZZ)})$ & $\lambda_F(p)=$ & $T(p)(\Delta_{19})$ & $\Tr(T(p))(\Delta^2_{25,19,5})$ \kern-0.4em \\& & $T(p)(\Delta_{25,5})$ & & \\\hline $2$ & $10176$ & $-96$ & $456$ & $6624$\\$3$ & $929988$ & $-527688$ & $50652$ & $90072$\\$5$ & $-36016170$ & $596139180$ & $-2377410$ & $-334979100$\\$7$ & $-40517568504$ & $-3608884496$ & $-16917544$ & $-31105966416$  \\\hline
\end{tabular}
\vskip10pt
\begin{tabular}{|c|c|c|c|}\hline $p$ & $\Tr(T(p)|V_{\mu'}^{\SO(9)(\ZZ)})$ & $T(p)(\Delta_{25,19,11,5})$ & $-T(p)(\Delta_{25,19,11,5})$\\ & & & $+(p^3+p^7)a_p(f)+\lambda_F(p)$\\\hline $2$ & $5280$ & $4800$ & $2^5.3^2.5.\mathbf{17}$\\ $3$ & $889920$ & $-302400$ & $-2^8.3^3.5.13.\mathbf{17}$\\$5$ & $-345413400$ & $-765121800$ & $2^{10}.3^2.5.\mathbf{17}.53.131$\\$7$ & $-29042227200$ & $29642547200$ & $2^9.3^3.5.7.\mathbf{17}.191.1459$\\\hline
\end{tabular}
\section{An endoscopic congruence for $\SO(9)$.}\label{SO9}
{\bf Example 6. }
Now consider a self-dual, cuspidal, automorphic representation of $\PGL_8(\A)$ discovered by Chenevier and Renard  \cite[Corollary**6.5, Table 8]{CR}, denoted $\Delta_{25,21,15,9}$ in their notation. By Arthur's symplectic-orthogonal alternative \cite[Theorem* 3.9]{CR}, it is the functorial lift of a discrete automorphic representation of $\SO(5,4)(\A)$. There is an associated $\tilde{\rho}:\Gal(\Qbar/\QQ)\rightarrow \GL_8(\QQ_q)$ (see \cite[Remark 7.6]{Sh}). The relevant space of algebraic modular forms for $\SO(9)$ is $3$-dimensional, spanned by Hecke eigenforms that are vectors in automorphic representations of $\SO(9)(\A)$ with Arthur parameters $\Delta_{25,21,15,9}$, $\Delta_{25}\oplus\Delta_{15}\oplus\Delta_{21,9}$ and $\Delta_{21,9}\oplus\Delta_{25,15}$ \cite[Table 1]{Me}. To get the Hecke eigenvalues we want, for $\Delta_{25,19,11,5}$, one must subtract the endoscopic contributions from the traces calculated by M\'egarban\'e for $p\leq 7$, in fact he has done that and listed the results in \cite[Table 7]{Me}. Here $F$ and $G$ are genus $2$ cuspidal, Hecke eigenforms for $\Sp_2(\ZZ)$, with $(j,k)=(8,8)$ and $(14,7)$ respectively.
\vskip10pt
\begin{tabular}{|c|c|c|c|}\hline $p$ & $T(p)(\Delta_{25,21,15,9})$ & $\lambda_F(p):=T(p)(\Delta_{21,9})$ & $\lambda_G(p):=T(p)(\Delta_{25,15})$ \\\hline $2$ & $-7200$ & $1344$ & $-3696$ \\
$3$ & $631200$ & $-6408$ & $511272$ \\ $5$ & $6175800$ & $-30774900$ & $118996620$ \\ $7$ & $25981995200$ & $451366384$ & $-82574511536$\\\hline
\end{tabular}
\vskip10pt
\begin{tabular}{|c|c|}\hline $p$ & $p^2\lambda_F(p)+\lambda_G(p)-T(p)(\Delta_{25,21,15,9})$\\
\hline $2$ & $2^4.3.5.\mathbf{37}$\\ $3$ & $-2^6.3.5^2.\mathbf{37}$\\ $5$ & $-2^8.3.5.\mathbf{37}.4621$\\ $7$ & $-2^7.3^3.5.\mathbf{37}.135197$\\\hline
\end{tabular}
\vskip10pt

The congruence verified in the second table for $p\leq 7$, if it
held for all $p$, would imply that (with $q=37$) $\,\,\tilde{\rhobar}$
has composition factors $\rhobar_F(-2)$ and $\rhobar_G$. The
irreducibility of $\rhobar_F$ and $\rhobar_G$ can presumably be
checked as in \cite[Proposition 4.10]{CL}. Now reasoning as in
Example~3, but in the opposite direction, we should expect to find
the prime factor $37$ in the numerator of $\pi^8\,\frac{L(1,\pi_F\otimes\pi_G)}{L(2,\pi_F\otimes\pi_G)}$, a ratio of critical
values for a degree-$16$ $\,\,\GSp_2\times\GSp_2$ $L$-function. (As
was pointed out to us by M. Watkins, $L(3,\pi_F\otimes\pi_G)$ is
not critical.)
Using the same techniques as in \S 4, we obtain the approximation
$\pi^8\, \frac{L(1,\pi_F\otimes\pi_G)}{L(2,\pi_F\otimes\pi_G)}\approx 6243.7501$,
likely to be correct to 7 or 8 significant figures.
This approximation is not as accurate as what we obtained in our
degree~8 examples, but we note that $6243.75 = \frac{3^3 \cdot 5^2 \cdot \mathbf{37}}{2^2}$.
This is suggestive of the expected prime factor of $37$.  However,
this evidence is not as convincing as, for example, equation~\eqref{eqn:manyzeros}.

Though it was the narrowness of the critical range that forced evaluation of the $L$-function at adjacent points, doing so had a fortuitous side-benefit. The quantity $\pi^8 \, \frac{L(1,\pi_F\otimes\pi_G)}{L(2,\pi_F\otimes\pi_G)}$ is small enough that 7 or 8 significant figures are enough to go comfortably beyond approximating just the integer part, and to reveal what is presumably the correct rational number. This would not have been the case for the much larger $\pi^{16}\, \frac{L(1,\pi_F\otimes\pi_G)}{L(3,\pi_F\otimes\pi_G)}$, even had
$L(3,\pi_F\otimes\pi_G)$ been critical.

\end{document}